\documentclass[11pt]{amsart}
\usepackage{amssymb}
\def\nat{{\bf N}}	
\def\C{{\bf C}}         
\def\Z{{\bf Z}}		
\def\set#1{\left\{#1\right\}}               
\def\abs#1{\left|#1\right|}		    
\def\tr{{\rm Tr}\,}                         
\def\rk{{\rm rk}\,}                         

\def\ferrbox{\rlap{\vrule height14pt}
\rlap{\vrule width14pt height0.4pt depth0pt}
\raise14pt\hbox{\vrule width14pt height0.4pt depth 0pt}
\rlap{\vrule height14pt}}
\newcount\ferrowcount
\def\ferrersrow#1{{\ferrowcount=0\hbox{$
\loop\ifnum\ferrowcount<#1 \advance\ferrowcount by1
\smash\ferrbox\vphantom{\vrule height14pt}\repeat$}}}

\def\cntr#1{{\setbox0=#1\lower0.5\ht0\copy0}}

\def\ferrpig#1{\ifx#1\end \let\next=\relax
 \else\ferrersrow{#1}\nointerlineskip\let\next=\ferrpig\fi\next}

\def\blackslug{\hbox{\hskip 1pt \vrule width 4pt height 8pt depth 1.5pt
	\hskip 1pt}}
\def\QED{\kern0.5em\blackslug\lower 8.5pt\null}
\def\qed{{\hfill\null\nobreak\hfill\hbox to 50pt{\quad\blackslug\lower 
8.5pt\null\hfil}}}
\def\jacqed{{\kern0.5em\lower3pt\hbox{$\Box$}}}
\newcounter{enumalpha}\renewcommand{\theenumalpha}{\alph{enumalpha}}
\newenvironment{enumal}{\begin{list}{(\theenumalpha)}{\usecounter{enumalpha}}}{\end{list}}

\def\ca{{\cal{A}}}

\def\cC{{\cal C}}

\theoremstyle{definition}
\newtheorem{definition}[equation]{Definition}
\newtheorem{example}[equation]{Example}
\theoremstyle{plain}
\newtheorem{proposition}[equation]{Proposition}
\newtheorem{lemma}[equation]{Lemma}
\newtheorem{theorem}[equation]{Theorem}
\newtheorem{corollary}[equation]{Corollary}

\def\cal{\mathcal} 
\def\Z{{\bf Z}} 
\input{xypic}
\xyoption{graph}
\newtheorem{rem}[equation]{Remark}
\def\md{{\cal D}}
\def\eqref#1{{$(\ref{#1})$}}
\def\nodiv{{\;/\hspace*{-0.4em}|\;}}
\def\bfcdot{{\hbox{$\scriptscriptstyle\bullet$}}}
\def\P{{\bf P}}

\def\roundup#1{{\lceil#1\rceil}}
\def\Reg{{\cal R}}
\title[equivariant euler characteristics for reflection groups]{Equivariant Euler characteristics of discriminants of reflection groups}
\author{Graham Denham}
\email{denham@noether.uoregon.edu}
\address{University of Oregon\\
Eugene, OR, USA 97403}
\thanks{Both authors acknowledge the support of 
NSERC Postdoctoral Fellowships} 
\author{Nicole Lemire}
\email{nlemire@uwo.ca}
\address{University of Western Ontario\\
London, Ontario, Canada N6A 5B7}
\date{\today}
\begin{document}
\begin{abstract}
Let $G$ be a finite, complex reflection group acting on a complex
vector space $V$, and $\delta$ its disciminant polynomial.
The fibres of $\delta$ admit commuting actions of $G$ and a cyclic
group.  The virtual $G\times C_m$ character given by
the Euler characteristic of a fibre is a refinement of the 
zeta function of the geometric monodromy, calculated 
in~\cite{delo95}.  We show that this virtual character is unchanged
by replacing $\delta$ by a slightly more general class of polynomials.
We compute it explicitly, by studying the poset of normalizers of
centralizers of regular elements in $G$, and the subspace arrangement
given by the proper eigenspaces of elements of $G$.  As a 
consequence we also compute orbifold Euler characteristics and
find some new ``case-free'' information about the discriminant.
\end{abstract}
\bibliographystyle{amsplain}
\maketitle
\section{summary}
Let $G$ be a finite reflection group acting on the vector space
$V=\C^{\ell}$.
Let $\ca$ denote the set of reflecting hyperplanes of $G$.  For each
$H\in\ca$, let $\alpha_H\in V^*$ be a linear functional with
kernel $H$.  The discriminant polynomial $\delta$ of $G$ is defined to be
$$
\delta=\prod_{H\in\ca}\alpha_H^{e_H},
$$
where $e_H$ is the order of the subgroup of $G$ that fixes $H$ pointwise.
$\delta$ is the $G$-invariant 
polynomial of smallest degree whose zero set is exactly
the set of reflecting hyperplanes.  Let $m=\deg\delta$.

The fibres of $\delta$ over $\C^*$ are diffeomorphic, by a theorem of
Milnor~\cite{milnorbook}; let $F=\delta^{-1}(1)$, the Milnor fibre of
$\delta$.  The action of $G$ on $V$ restricts to an action on $F$.  
At the same time, a cyclic group $C_m$ acts on $F$, generated by a
geometric monodromy map $h:F\rightarrow F$ defined by $h(x)=e^{2\pi i/m}x$.

The actions of $G$ and $C_m$ commute.  Let $\Gamma=G\times C_m$.  Then
$H_\bfcdot(F,\C)$ is a finite-dimensional representation of $\Gamma$.
In this paper, we consider the Euler characteristic of $F$, valued in
the character ring of $\Gamma$.  That is, we define a virtual character
$\chi_\Gamma$ by
\begin{equation}\label{eq:ecdef}
\chi_\Gamma(F)(g)=\sum_{p\geq0}(-1)^p\tr(g,H_p(F,\C)).
\end{equation}

This is a refinement of the usual zeta function of the monodromy of 
the Milnor fibre, which Denef and Loeser~\cite{delo95} have calculated for 
all reflection groups.  Their technique uses Springer's theory of regular
elements~\cite{spr74}: $g\in G$ is called a {\em regular element} of $G$ iff
it has an eigenvector that is not contained in any reflecting hyperplane.

In particular, Springer~\cite{spr74} has shown
that the centralizer of a (noncentral) regular element in $G$ acts
as a reflection group on a (proper) subspace of $V$.  Using this
idea and its elaboration in~\cite{delo95,lespr99}, we find a recursive
formula for the Euler characteristic $\chi_\Gamma$ (Theorem~\ref{th:main}).

For a fixed reflection group $G$, let 
$M_G=M=V-\cup_{H\in \ca}H$ be the hyperplane complement, and $U$ its image
in $\P(V)$.
These spaces have been studied extensively in the context of hyperplane
arrangement theory; see, for example, \cite{leh95}.  On the other hand,
let ${\cal E}$ denote the set of all maximal eigenspaces $E$ of elements
of $G$ for which $E\subsetneq V$.  Let
\begin{equation}\label{eq:mcirc}
M^\circ=V-\bigcup_{E\in {\cal E}}E,
\end{equation}
and $U^\circ$ the image of $M^\circ$ in $\P(V)$.
$U^\circ$ is the complement of a projective subspace arrangement, 
in the sense of Bj\"orner~\cite{bjo94}, and to the authors' knowledge
has not been studied directly before.  

$G/Z(G)$ acts freely on $U^\circ$.  We find a formula for
the Euler characteristic
of the orbit space in terms of degrees, codegrees,
and regular numbers \eqref{eq:inverse}, and we calculate it
for each irreducible $G$ (Theorem~\ref{th:cases}).  
We show that this determines 
$\chi_\Gamma$ for each $G$ (Theorem~\ref{th:main}).  

\section{Springer's Theory of regular elements}\label{sec:spr}
In this section, we recall the theory of regular elements
and set up our notation.  We refer to~\cite{ot} for background
on reflection groups and hyperplane arrangements and to~\cite{spr74}
for background on the theory of regular elements.

Let $V=\C^\ell$ and let $G$ be a
finite reflection group acting on $V$. We will denote by 
$\C[V]$ the algebra of polynomial functions on $V$.
The {\em degrees} $d_1,\ldots,d_\ell$ of $G$ are the degrees of any set of homogeneous polynomials
which generate the $G$-invariant polynomial ring $\C[V]^G$.
The order of $G$ and of its centre $Z(G)$ are determined
in terms of its degrees:
\begin{equation}
|G|=\prod_{i=1}^{\ell}d_i\qquad |Z(G)|=\gcd\set{d_i}
\end{equation}

A vector $v\in V$ is called {\em regular} if it is 
not contained in a reflection hyperplane of $G$.  
An element $g\in G$ is called
{\em regular} if it has a regular eigenvector.  
Let $g\in G$ be regular of order
$d$. Let $v$ be a regular eigenvector with corresponding eigenvalue
$\xi$ and let $V(g,\xi)$ denote the $\xi$-eigenspace of $g$.
We will refer to $(g,\xi)$ as a {\em regular} ($d$-){\em pair}.

With this notation, we have:
\begin{theorem}[Springer~\cite{spr74}]\label{th:sprth}~
\begin{enumerate}
\item[(a)] The root of unity $\xi$ has order $d$.
\item[(b)] $V(g,\xi)$ has dimension $a(d)=\abs{\set{i:d|d_i}}$.
\item[(c)] The centralizer $C_G(g)$ is a reflection group in $V(g,\xi)$ 
whose degrees are $\set{d_i:d|d_i}$ and whose order is $\prod_{d|d_i}d_i$.
\end{enumerate}
\end{theorem}

The orders of the regular elements of $G$ are called the {\em regular numbers}
of $G$.  Let $\Reg$ denote the poset of regular numbers, ordered by
divisibility.

The group $G\le GL(V)$ also acts naturally on the algebra of polynomial 
vector fields on $V$, $\C[V]\otimes V$.  
The module $(\C[V]\otimes V)^G$ is free over $\C[V]^G$.  
Following \cite{bmr98}, 
the {\em codegrees} $d^*_1,\ldots,d^*_\ell$ are defined to be 
the degrees of a homogeneous basis, with the convention that derivations 
have degree $-1$.
By a theorem of Orlik and Solomon~\cite{os84}
\begin{equation}\label{eq:poin}
\sum_{i=1}^l\dim H^i(U,\C)t^i=\prod_{i=2}^\ell(1+(d^*_i+1)t)
\end{equation}

Using a case-based argument,
Denef and Loeser~\cite[Theorem 2.8]{delo95} proved that, for a 
regular $d$-pair $(g,\xi)$, the codegrees of $C_G(g)$ acting on $V(g,\xi)$
are 
\begin{equation}\label{eq:codeg}
\set{d^*_i:d|d^*_i}
\end{equation}  Lehrer and Springer~\cite[Theorem~C]{lespr99}
later reproved this result in a case-free way.

\section{euler characteristics}\label{sec:ec}
Following \cite{bmr98}, for each $G$-orbit of hyperplanes $\cC\in\ca/G$,
set
$$
\delta_\cC=\prod_{H\in\cC}\alpha_H^{e_\cC},
$$
where $e_\cC$ is defined to be the common value of $e_H$ for all $H\in\cC$.
noting that $e_H$ is constant for all $H\in\cC$.
Consider any homogeneous,
$G$-invariant polynomial $f\in\C[V]^G$ with zero locus equal to
$\bigcup_{H\in\ca}H$.  Then $f$ has the form
\begin{equation}\label{eq:f}
f=\prod_{\cC\in\ca/G}\delta_\cC^{a_{\cC}}
\end{equation}
for some positive integers $a_{\cC}$; in particular, the discriminant is
obtained by choosing all $a_{\cC}=1$.  We shall call such $f$ 
{\em unreduced discriminant polynomials}.

Denote the degree of $f$ as above by $m$.  Then $f(gv)=f(\zeta v)=\zeta^m f(v)$
for any regular element $g$ with eigenvalue $v$, so the order of 
$\zeta$ must divide $m$.  That is, all regular numbers $d\in\Reg$
divide $m$.

Let $F=f^{-1}(1)$.  Let $P=\pi_1(M,1)$, the pure Artin braid group 
corresponding to the group $G$.  Since $F$ is homotopy-equivalent to 
an infinite cyclic cover of $M$, we have
$H_\bfcdot(F,\C)\cong H_\bfcdot(M,\C[t,t^{-1}])$, where 
$\C[t,t^{-1}]\cong \C\uparrow_{\pi_1(F)}^P$ as a $P$-module; see 
\cite[Section~2.1]{d97} and \cite{de01}.
Explicitly, $P$ has a set of generators $\set{\gamma_H:H\in\ca}$
for which $\gamma_H$ acts by multiplication by $t^{a_\cC e_\cC}$.

The complement $M$ 
is known to be a $K(P,1)$ space for all irreducible reflection
groups with the possible exception of $G_{24}$, $G_{27}$, $G_{29}$, $G_{31}$,
$G_{33}$, and $G_{34}$; see \cite[(2.11)]{bmr98} for references.  Thus
we have $H_\bfcdot(F,\C)\cong H_\bfcdot(P,\C[t,t^{-1}])$ except 
perhaps in these cases.  Let $B=\pi_1(M/G,1)$, the braid group.
Then $H_\bfcdot(P,\C[t^{\pm1}])^G\cong H_\bfcdot(B,\C[t^{\pm1}])$ (over
$\C$).  For all real reflection groups, this is computed 
explicitly in \cite{cpss99}.

The Lefschetz zeta function of $F$ is defined to be
$$
Z(F)=\prod_{p\geq0}\det(1-h_*t|_{H_p(F,\C)})^{(-1)^p},
$$
where $h^*$ denotes a preferred generator of $C_m$ acting in homology.
Since a complex representation of $C_m$ is determined by the
characteristic polynomial of a generator of the group, the zeta function
can be seen as the restriction to $C_m$ of the Euler characteristic
$\chi_\Gamma(F)$ defined in \eqref{eq:ecdef}, written multiplicatively.
We will identify $C_m=\langle h^*\rangle$ with the cyclic group of 
$m$ elements in $\C^*$.  For convenience, we will take the convention
that $\alpha\in C_m$ acts on $F$ by multiplication by $\alpha^{-1}$.

For a reflection group $G$ and integers $d|m$, define
$$
I_d(G)=I_d=1\uparrow_{C_d}^{G\times C_m}
$$
whenever $C_d$ is a cyclic subgroup of order $d$ generated by a regular
pair $(g,\zeta)$ of order $d$. 

\begin{lemma}\label{lem:id}
The cyclic groups generated by any two regular pairs of the same order
are conjugate in $\Gamma$.  In particular,  
the definition of $I_d(G)$ above does not depend on the choice of 
regular pair $(g,\zeta)$.  
The normalizer in $\Gamma$ of the cyclic subgroup generated
by $(g,\xi)\in \Gamma$ is $C_G(g)\times C_m$.
\end{lemma}
\begin{proof}
Let $(g,\zeta)$ and $(g',\xi)$ be two 
regular pairs of 
order $d$, generating cyclic groups $K$ and $K'$, respectively.  Since any 
two primitive
roots of unity of the same order generate the same cyclic subgroup of $\C^*$,
$\xi=\zeta^k$ for some $k$.  Then $(g^k,\xi)$ is also a regular pair of 
order $d$, so $g^k$ and $g'$ are conjugate, by \cite[4.2]{spr74}.
It follows that the subgroups $K$ and $K'$ are conjugate.
$I_d(G)$ is well defined since the permutation characters induced
from $K$ and $K'$ are the same~\cite[10.12]{CRvolI}.
Since
$(h,\alpha)\in N_{\Gamma}\langle(g,\xi)\rangle$  iff 
$(hgh^{-1},\alpha)=(g^k,\xi^k)$ for some $k$ iff  $h\in C_G(g)$, 
$\alpha\in C_m$,
we have that $N_{\Gamma}(\langle (g,\xi)\rangle)=C_G(g)\times C_m$.
\end{proof}

The following theorem appeared independently as 
\cite[Theorem 2.5]{delo95} and \cite[Corollary 5.8]{leh95}.
\begin{lemma}\label{lem:czrs}
If $(g,\zeta)$ is a regular pair and $V$ is
the $\zeta$-eigenspace of $g$, then the centralizer $C(g)$ acts as
a reflection group on $V$.  Its reflecting hyperplanes are $\set{H\cap V:
H\in\ca}$.
\end{lemma}

\begin{definition}
For a regular $d$-pair $(g,\zeta)$, we define 
$U(g,\zeta)$ as the projective hyperplane complement for $C_G(g)$
acting on $V(g,\zeta)$ where
$(g,\zeta)$ is a regular $d$-pair.

From Lemma~\ref{lem:id},
$C_G(g)$ is conjugate to $C_G(g')$ if $(g,\zeta)$ and $(g',\xi)$ are 
both regular $d$-pairs.  This means that they are isomorphic as reflection
groups, and we will refer to them, up to isomorphism, as $G(d)$,
as in~\cite{lespr99}.

By the lemma above, then, $U(g,\zeta)$ and $U(g',\xi)$ are diffeomorphic;
we will refer to them as $U_G(d)$.
\end{definition}

\begin{definition}\label{df:md}
Define a poset $\md=\md_G$ by
$$
\md=\set{d:d=\abs{Z(C_G(g))}\hbox{~for a regular element $g$}},
$$
ordered by divisibility.  
That is, $\md$ is the set of orders of regular elements
$g$ that are maximal with respect to the property of having a given 
centralizer.  For elements $d\in\Reg_G$, define $\roundup{d}$
to be the least multiple of $d$ in $\md$.
\end{definition}

Note that $\set{G(d): d\in\md}$ forms a complete set of representatives
of the isomorphism classes of centralizers of regular elements.

Recall $M^\circ\subseteq M$ from~\eqref{eq:mcirc}, and $U^\circ\subseteq U$.  
By construction,
\begin{proposition}\label{prop:freely}
$G/Z(G)$ acts freely on $U^\circ$.
\end{proposition}

We can now state our main result.

\begin{theorem}\label{th:main}
Let $G$ be a reflection group, $f$ an unreduced discriminant polynomial
\eqref{eq:f} of degree $m$, and $F$ its Milnor fibre.
Then
$$
\chi_\Gamma(F)=\sum_{d\in\md}a_d I_d,
$$
where the integers $a_d$ are given by
\begin{equation}\label{eq:defofa}
a_d=\chi(U(d)^\circ/G(d))
\end{equation}
\end{theorem}

A case analysis gives a more refined description.  First, 
$\chi_\Gamma(F)$ is zero unless $G$ is irreducible, since $\chi(U)=0$
in this case: this appears first in the language of matroids in \cite{crapo}.
Denef and Loeser~\cite[2.9]{delo95} show that the centralizers of 
regular elements in irreducible $G$ are themselves irreducible.  With
this in mind, it is enough to calculate $a_z=\chi(U^\circ/G)$, 
where $z=\abs{Z(G)}$, for each irreducible $G$.  We obtain:

\begin{theorem}\label{th:cases}
For an irreducible reflection group $G$ of rank $n$, 
$\chi(U^\circ/G)=(-1)^{n-1}$ if
$G$ is in the list below.  Otherwise, $\chi(U^\circ/G)=0$.
\begin{enumal}
\item Irreducibles of rank $\leq2$;
\item Irreducibles of the form $G(de\ell,e\ell,\ell)$, except $G(3,3,3)$;
\item $G(de\ell,e\ell,2\ell)$ where $e\ell$ is odd;
\item Exceptionals $G_{29}$ and $G_{34}$.
\end{enumal}
\end{theorem}

\begin{corollary}
For a Milnor fibre of a given reflection group as above,
$\chi_\Gamma$ is a linear combination of at most six permutation characters
$I_d$, with coefficients $\pm1$.
\end{corollary}

The value of $\chi_\Gamma$ for each irreducible reflection group is
tabulated in Section~\ref{sec:calc}.

We also observe empirically that, like the zeta function, $\chi_\Gamma$
continues to be a braid diagram invariant, in the sense of Brou\'e,
Malle, and Rouquier~\cite{bmr98}.

\begin{example}\label{ex:e8}
Let $G$ be the irreducible reflection group of type $E_8$.  The degrees
are $[2,8,12,14,18,20,24,30]$ and
the poset $\md$ is:
$$
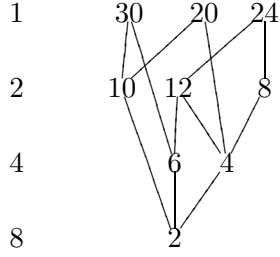
\xygraph{[]~{<1cm,0cm>:<0cm,0cm>}
*{30}="A"[]!{!(1,0)}
*{20}="B"[]!{!(0.8,0)}
*{24}="C"-@{-}[]!{!(0,-1)}
*{8}="D"
"C"-@{-}[]!{!(-1.15,-1)}
*{12}="E"
"A"-@{-}[]!{(-.1,-1)}
*{10}="F"-@{-}"B"
"E"-@{-}[]!{!(-0.05,-1)}
*{6}="G"-@{-}"A"
"E"-@{-}[]!{!(0.65,-1)}
*{4}="H"-@{-}"D"
"H"-@{-}"B"
"H"-@{-}[]!{!(-0.7,-1)}
*{2}="I"-@{-}"G"
"I"-@{-}"F"
"A"[]!{!(-1.5,0)}
*{1}[d]*{2}[d]*{4}[d]{8}
}
$$
By the Shephard-Todd classification, the centralizers of elements of
order $4$ and $6$ are, respectively, $G_{31}$ and $G_{32}$.  Applying
Theorems~\ref{th:main} and \ref{th:cases} shows that
$$
\chi_{E_8\times C_m}=I_{30}+I_{24}+I_{20}-I_{12}-I_{10}-I_8.
$$
\end{example}

The rest of this paper is as follows.
Theorem~\ref{th:main} is proven in Section~\ref{sec:main}.  
Theorem~\ref{th:cases}
is proven in Section~\ref{sec:calc}.  We are unable to provide a 
case-free proof of Theorem~\ref{th:cases} in general, although in 
Section~\ref{sec:rank2} we do so for rank $2$ irreducibles.

\section{Proof of Theorem~\ref{th:main}}\label{sec:main}
Throughout this section, fix a particular reflection group $G$
acting on $V=\C^\ell$.  We will use the notation of Section~\ref{sec:ec}.

For a regular element $g\in G$ of order $d$, let $i_G(d)$ be the order of
its centralizer, and let $u_G(d)$ be the Euler characteristic of 
$U(d)$.  Then $i(d)=\prod_{d|d_i}d_i$ (Theorem~\ref{th:sprth}).
By evaluating~\eqref{eq:poin} at $t=-1$ and using~\eqref{eq:codeg}, 
we see $u(d)=\prod_{d|d^*_i,i>1}(-d^*_i)$.

If $(g,\zeta)$ is a regular pair of order $d$, then by \cite[3.4]{spr74},
there is a regular pair $(h,\zeta)$ of order $\roundup{d}$ for which 
$C(g)=C(h)$ (and, hence, $U(g,\zeta)=U(h,\zeta)$).  It follows that
$i(d)=(\roundup{d})$ and $u(d)=u(\roundup{d})$.

For $(g,\zeta)\in G\times C_m$, it follows from a refinement
of the Hopf-Lefschetz fixed
point theorem~\cite{kbrown82} that 
\begin{eqnarray}
\chi_\Gamma(F)(g,\zeta)&=&\chi(F^{(g,\zeta)})=\deg(\delta)\cdot\chi(U(g,\zeta))\nonumber\\
&=&\begin{cases}m\cdot u(\roundup{d})&\text{if $(g,\zeta)$ is regular 
of order $d$;}\\
0&\text{otherwise.}\\ \end{cases}\label{eq:cg}
\end{eqnarray}
\begin{lemma}\label{lem:exists}
There exist rationals $\set{a_d:d\in\md}$ for which 
\begin{equation}\label{eq:sumofids}
\chi_\Gamma=\sum_{d\in\md}a_d I_d.
\end{equation}
\end{lemma}
\begin{proof}
We claim that, if $I_d$ is induced from a regular pair $(g,\zeta)$ 
of order $d$, then
\begin{equation}\label{eq:exists}
I_d(h,\xi)=\begin{cases}
\frac{m}{d}\abs{C(h)}&\text{if $(h,\xi)$ is regular of order $d'$, $d'|d$;}\\
0&\text{otherwise.}\end{cases}
\end{equation}
This follows by directly evaluating the induced character and using
Lemma~\ref{lem:id}: $I_d(h,\xi)$ is
nonzero only if $(h,\xi)$ is conjugate to a power of $(g,\zeta)$, the generator
of $C_d$, which is equivalent to $(h,\xi)$ being
a regular pair of an order dividing that of $g$.  

Define an equivalence relation $\sim$ on $G\times C_m$ by
setting $(g,\zeta)\sim(h,\xi)$ iff either: neither is a regular pair; or,
both are regular pairs, and $\roundup{g}=\roundup{h}$.

The equivalence classes of $\sim$ are unions of conjugacy
classes.  Moreover, \eqref{eq:exists} shows that each $I_d$ is
constant on classes of $\sim$.  

The characters $I_d$ span the ${\bf Q}$-vector space of functions that
are both constant on $\sim$ classes and zero on the nonregular equivalence
class, since their values on the regular pairs form a triangular matrix,
by \eqref{eq:exists}.

By \eqref{eq:cg}, $\chi_\Gamma$ is such a class function, which completes
the proof.
\end{proof}

\subsection{coefficients $a_d$ from (co)degrees}
By evaluating \eqref{eq:sumofids} on a regular pair of order $d\in\md$,
we obtain
$$
m\cdot u(d)=\sum_{k\in\md:d|k}a_k\frac{m}{k}i(d),
$$
whence by M\"obius inversion,
\begin{eqnarray}
a_d&=&d\sum_{d|k}\mu(d,k)\frac{u(k)}{i(k)},\label{eq:inverse}\\
&=&d\sum_{d|k}\mu(d,k)\prod_{i:d|d^*_i\atop d^*_i\neq0}
d^*_i\prod_{i:d|d_i}d_i^{-1}\nonumber
\end{eqnarray}
where $\mu$ is the M\"obius function of the poset $\md$.

\begin{rem}
Any rational character of a finite group can be
expressed as a rational linear combination of induced permutation characters
from cyclic subgroups by the Artin-Brauer 
induction theorem~\cite[15.2]{CRvolI}.  The coefficients
can be determined and are non-zero only if the character is non-zero on that
cyclic subgroup.  In our situation, however, the cyclic subgroups $C_d, d\in \md$ are representatives of isomorphism classes (hence a subset of the
set) of cyclic subgroups 
of $\Gamma$
with $\chi_{\Gamma}(F)$ non-zero.  Our argument above is then a
direct way to compute the coefficients $a_d,d\in \md$ for our special case.
\end{rem}

\subsection{induction from a centralizer subgroup}
Given a regular element $g_0\in G$ of order $e\in\md$, let $H=C_G(g_0)$.
From \cite[4.2]{spr74}, and Definition~\ref{df:md} we observe:
\begin{lemma}
The maximal regular numbers of $H$ are
$\md_H=\set{d\in\md_G:e|d}$.
\end{lemma}

By Lemma~\ref{lem:exists}, 
\begin{eqnarray*}
\chi_{H\times C_m}&=&\sum_{d\in\md_H}a'_d I_d(H),\quad\hbox{and}\\
\chi_{G\times C_m}&=&\sum_{d\in\md}a_d I_d
\end{eqnarray*}
for some coefficients $\set{a'_d}$, and $\set{a_d}$.

\begin{theorem}\label{th:compatible}
If $G$ and $H$ are as above, and $m'|m$, then
$$
\chi_{G\times C_m}=\chi_{H\times C_m'}\uparrow_{H\times C_{m'}}^{G\times C_m}
+\sum_{d\in\md:e\nodiv d}a_d I_d.
$$
Consequently, 
in the notation above, $a^{\phantom{1}}_d=a'_d$ for all $d\in\md_H$.
\end{theorem}
\begin{proof}[Proof of theorem]
We need to show that the values $a_d$ given by equation \eqref{eq:inverse} 
for $d\in\md_H$
are the same in $H$ as they are in $G$.  
Specifically, we need to show for multiples $k$ of $e$, that
$i_G(k)=i_H(k)$, and that $u_G(k)=u_H(k)$.

So suppose that $g\in H$ is a regular
element of order $k$, where $e|k$.  
Then $g_0$ is conjugate to a power of $g$,
so $C_G(g)\subseteq C_G(g_0)=H$; that is, $C_G(g)=C_H(g)$, so $i_G(k)=i_H(k)$.

Now assume, without loss of generality, that $g=g_0^r$ for some $r$, and
choose $\zeta$ so that $(g,\zeta)$ is a regular pair.  
Then the $\zeta$ eigenspace of $g$ in $V_G$ is contained in $V_H$, so 
$U_G(g,\zeta)=
U_H(g,\zeta)$ (Lemma~\ref{lem:czrs}.)  Then $u_G(k)=u_H(k)$ as well.
\end{proof}
\subsection{interpretation of coefficients $a_d$}
For $d\in\md$, let $F_d=F\cap V_d$, where 
$$
V_d=\bigcap_{d\nodiv d_i}f_i^{-1}(0),
$$
the variety of eigenvectors of elements of $g$ having eigenvalue
a primitive $d$th root of unity; see~\cite[3.2]{spr74}.  
In particular, $F_z=F$, where $z$ is the order of the centre of $G$.
Let 
$$
F_d^\circ=F_d-\bigcup_{d|d'\atop d'\neq d}F_{d'}.
$$
The following lemma is standard.
\begin{lemma}\label{lem:ecs}
Suppose $X\subseteq M$ is a smooth embedding of complex manifolds without
boundary and $X$ is closed in $M$.  Putting $U=M-X$, we have
$$
\chi(X)+\chi(U)=\chi(M).
$$
Furthermore, if a finite group $G$ acts on $M$ so that $GX=X$, we have
$$
\chi_G(X)+\chi_G(U)=\chi_G(M).
$$
\end{lemma}
\begin{proposition}\label{prop:additive}
For any $d\in\md$,
$$
\chi_{G\times C_m}(F_d)=\sum_{d'\in\md\atop d|d'}\chi_{G\times C_m}
(F_{d'}^\circ).
$$
\end{proposition}
\begin{proof}
Use induction on the
dimension of $F_d$.  If minimal, $F_d=F_d^\circ$ and there is nothing to
prove.  Otherwise, we have
\begin{eqnarray*}
F_d&=&F_d^\circ\cup\bigcup_{d|d'\atop d'\neq d}F_{d'}^\circ\\
&=&F_d^\circ\cup\bigcup_{d'\succ d}F_{d'}
\end{eqnarray*}
where ``$\succ$'' denotes the covering relation in the poset $\md$,
and all of the unions are disjoint.
The proof is completed by 
applying Lemma~\ref{lem:ecs} and the induction hypothesis to the second
equality.
\end{proof}

To complete the proof of Theorem~\ref{th:main}, it remains to show that
the coefficients $a_d$ given by Lemma~\ref{lem:exists} satisfy
$a_d=\chi(U(g,\zeta)^\circ/C_G(g))$ for regular pairs $(g,\zeta)$ of
order $d\in\md$.  

Let $z=\abs{Z(G)}$.
From Proposition~\ref{prop:freely}, the quotient map 
$F_z^\circ\rightarrow U^\circ/G$ is a covering with deck transformation
group $G\times C_m/\langle(g,\zeta)\rangle$, where $(g,\zeta)$ is any
regular $z$-pair.  It follows that 
\begin{eqnarray*}
\chi_{G\times C_m}(F^\circ)&=&\chi(U^\circ/G)\cdot
1\uparrow_{C_z}^{G\times C_m},\\
&=&\chi(U^\circ/G)\cdot I_z.
\end{eqnarray*}

By Theorem~\ref{th:compatible}, then, for each $d\in\md$,
$$
\chi_{\Gamma}(F^\circ_d)=\chi(U(d)^\circ/G(d))I_d,
$$

Now, using Proposition~\ref{prop:additive} with $d=z$, we have 
\begin{eqnarray*}
\chi_{\Gamma}(F)
&=&\sum_{d\in\md}\chi_{\Gamma}(F^\circ_d)\\
&=&\sum_{d\in\md}a_d I_d;
\end{eqnarray*}
equating the coefficient of $I_d$ for each $d$ gives the characterization
of the values $a_d$ that we claimed.


\section{The rank two case}\label{sec:rank2}
For finite groups acting on $\C^2$, a stronger version of 
Theorem~\ref{th:cases}(a) is obtained from the theory of du Val
singularities, for which we refer to \cite{dur79}.
\begin{theorem}
Let $G$ be an irreducible finite subgroup of $U_2(\C)$. 
Then 
$$\chi(U_G^\circ/G)=-1.$$
\end{theorem}

\begin{proof}\def\oH{\overline{H}}
Note that
a finite subgroup $G$ of $GL_\ell(\C)$ can always be embedded in $U_\ell(\C)$.
Let $Z:=Z(U_\ell(\C))=\set{\alpha I: \alpha\overline{\alpha}=1}$.
Let $H:=(G\cdot Z)\cap SU_\ell(\C)$; then
$Z(H)=H\cap SU_\ell(\C)$.
It is easy to check that the map $G/Z(G)\to H/Z(H)$
defined by $gZ(G)\mapsto g(\det(g)^{1/\ell})Z(H)$
is an isomorphism; therefore
$G$ and $H$ are both central extensions of the same group $\oH$.
Then $G$ and $H$ act on $\C^\ell$ and $\oH$ acts on
$\P^{\ell-1}$.  Proper eigenspaces for the action of $G$ on $\C^\ell$
are the same as those for $H$.  Moreover,
if $v$ is an eigenvector for an element $h\in H$ , its image
$[v]\in \P^{\ell-1}$ is a fixed point of $hZ(H)\in \oH$.

Now, let $\ell=2$.  Note that for $G\le U_2(\C)$, there
exists a reflection group $G'\le U_2(\C)$ with $\oH= G/Z(G)\cong 
G'/Z(G')$.  Then $\P^1/\oH\cong\P^1/G\cong \P^1/G'\cong \P^1$ where the last
isomorphism follows from the Shephard Todd Chevalley theorem.

Since $PSL_2(\C)\cong SO_3$, $\oH$ is isomorphic
to a finite subgroup of $SO_3$.
The finite subgroups of $SO_3$ are known:
they are the groups of symmetry of the regular polyhedra:
cyclic, dihedral, tetrahedral ($A_4$), octahedral ($S_4$) and icosahedral
($A_5$).  Since $G$ was assumed to act irreducibly, $\oH$ is not cyclic.
For the remaining finite subgroups of $SO_3$, Klein~\cite{klein} showed
that there are exactly three orbits of points in $S^2$ of points with
nontrivial stabilizers (corresponding to vertices, barycenters of 
edges of the regular polyhedron, and barycenters of faces.)
Since $\P^1/G\cong\P^1$, it follows that
$U^\circ_G/G=\P^1-\set{p_0,p_1,p_2}$, where the $p_i$'s are 
the 3 ``special'' orbits under the action of $H$ on $\P^1$.
Thus $\chi(U^\circ_G/G)=2-3=-1$.
\end{proof} 

\section{Proof of Theorem~\ref{th:cases}}\label{sec:calc}
The Shephard-Todd classification of irreducible
(complex) reflection groups consists of one infinite family
and 34 exceptional groups labelled as $G_4,\ldots, G_{37}$.
The tables in Figures~\ref{figa},\ref{figb},\ref{figc}
give $\chi_\Gamma$ for the exceptional
groups; these values are readily calculated from \eqref{eq:inverse}.

For any reflection group $G$, put $c(G)=\chi(U^\circ/G)$.  
From Theorem~\ref{th:main},
this equals the coefficient of $I_{|Z(G)|}$ in $\chi_\Gamma$.
By Theorem~\ref{th:compatible}, for $d\in \md$,
the coefficient $a_d$ of $I_d$ is $c(G(d))$.  

Theorem~\ref{th:cases} claims that for exceptional irreducible $G$,
$c(G)=0$ unless $G$ has rank $\ell=2$ or $G$ is one of $G_{29}$ or $G_{34}$,
in which case $c(G)=(-1)^{\ell-1}$.  Our proof in ranks $\ell>2$ is
by inspection, after having computed $\chi_\Gamma$ for each group.

The rest of this section is devoted to proving Theorem~\ref{th:cases}
for the infinite family of reflection groups $G(r,p,\ell)$.

\subsection{(co)degrees and regular numbers}
For $r,p\in\nat$ with
$p|r$ and rank $\ell\geq 2$, $G(r,p,\ell)$ is a group of 
order $r^\ell/p\ell!$.  It is the semidirect
product of the symmetric group $S_\ell$ acting by permutations on the
standard basis $\set{e_i:1\leq i\leq \ell}$, and the group of diagonal maps 
$e_i\mapsto \theta_ie_i$, where $\theta_i^r=1$ and
$(\theta_1\cdots\theta_\ell)^q=1$ where $q=r/p$.  
The group acts irreducibly on $\C^\ell$ iff $r>1$ and 
$(r,p,\ell)\ne (2,2,2)$.  
$G(r,1,\ell)$ is the 
full monomial group $C_r^\ell\rtimes S_\ell$ and that the Weyl groups 
$A_{\ell-1},B_\ell, D_\ell, G_2$ and the dihedral groups $I_2(\ell)$ 
equal respectively
$G(1,1,\ell), G(2,1,\ell), G(2,2,\ell), G(6,6,2)$ and $G(\ell,\ell,2)$.
The codegrees are calculated in \cite{os80b}, and the regular numbers appear
in \cite{co76}.  
The degrees of $G(r,p,\ell)$ are 
\begin{equation}\label{eq:degs}
\begin{cases}
r,2r,\ldots,(\ell-1)r,\ell q, &p|r, r>1\\
2,3\ldots,(\ell-1),\ell, &p=r=1.
\end{cases}
\end{equation}
the order of the center $z$ is $q\gcd(r,\ell)$.  The codegrees are
\begin{equation}\label{eq:codegs}
\begin{cases}
0,r,2r,\ldots,(\ell-1)r, & p<r\\
0,r,2r,\ldots,(\ell-2)r,(\ell-1)r-\ell & p=r> 1, \ell>1\\
0,1,2 ,\ldots,(\ell-2), &p=r=1,\ell>2
\end{cases}
\end{equation}
Note that the adjustments in the degrees and codegrees
of  $G(1,1,\ell)=S_\ell$ are made so that $S_\ell$ acts irreducibly.

\begin{rem}
Since the degrees and codegrees of $G(d)$ are those of $G$ which are 
divisible by $d$, we see for $G=G(r,p,\ell)$ and $e=(\gcd(d,r))^{-1}d$,
the degrees of $G(d)$ are
\begin{equation}
\begin{cases}\label{gddegs}
er,2er,\ldots,\lfloor\frac{\ell-1}{e}\rfloor er, & d\nodiv\ell q \\
er,2er,\ldots,\lfloor\frac{\ell-1}{e}\rfloor er,\ell q, &d|lq\\
2,3,\ldots,\ell, &d=r=1,
\end{cases}
\end{equation}
as noted in~\cite{lespr99}.
Lehrer and Springer prove indirectly in~\cite[5.2]{lespr99} that $G(d)=G(r',p',\ell')$ for
$G=G(r,p,\ell)$ and $d\in\Reg$.  In the proposition below we make the
determination of $G(d)$ explicit to help in our computation of $c(G)$.
\end{rem}

\begin{proposition}\label{prop:regnos}
\noindent
\begin{enumerate}
\item[(a)] For $G=G(r,p,\ell)$, $q=r/p>1$, we have $\Reg=\set{d:d|\ell q}$
and $\md=\set{kq: t|k|\ell}$ where $t=\gcd(p,\ell)$. Then 
for $kq\in \md$, $G(kq)=G(kr/t,p,\ell t/k)$.
\item[(b)] For $G=G(r,r,\ell)$, $\ell>1$, we have $\Reg=\set{d:d|\ell q} 
\cup\set{d:d|(\ell-1)r}$ and 
$$
\md= \set{kr:  k|\ell-1} \cup \set{d: z \vert d\vert \ell}
$$
where $z=\gcd(r,\ell)$.
For $k|\ell-1, kr\ne z$, we have $G(kr)=G(kr,1,(\ell-1)/k)$
and for $z|d|\ell$, we have $G(d)=G(dr/z,r,\ell z/d)$.
\end{enumerate}
\end{proposition}

\begin{proof}
In each case, let $z=\abs{Z(G)}$, equal to the greatest common divisor of the
degrees.  In the first case, $z=\gcd(pq,\ell q)=tq$, where $t=\gcd(p,\ell)$.
In the second case, this is just $z=\gcd(r,\ell)$.

\noindent Case (a):
For $G=G(r,p,\ell)$ with $q>1$, the only maximal regular degree
of $G$ is $\ell q$~\cite[2.11]{co76}.  
Thus the set of regular numbers
are $\Reg=\set{d:d|\ell q}$.
If $d\in \md$, then $d$ is a $\gcd$ of a subset of the degrees.
So $qt|d$, since $z=qt$ is the gcd of all the degrees.
To show that 
$$
\md=\set{kq: \gcd(p,\ell)|k|\ell}
$$
we need to show that $G(d)$ are distinct for distinct $d$ in this set.
Observe that if $\gcd(p,\ell)|k|\ell$, then $\gcd(p,k)=\gcd(p,\ell)$.
But as is implicitly shown in the proof of \cite[5.2]{lespr99}
we have for $t|k|\ell$,
$$G(kq)=G(kr/t,p,\ell t/k)$$ 
So we have shown that $\md$ and $G(d), d\in \md$ are described as in the claim.

\noindent Case (b):
For $G=G(r,r,\ell)$, with $\ell>1$, 
the maximal regular degrees are $(\ell-1)r$ and $\ell$~\cite[2.11]{co76}.
So the set of regular numbers is 
$$
\Reg=\set{d: d|\ell}\cup \set{d: d|(\ell-1)r}
$$
For $d\in \md$, $z=\gcd(r,\ell)$, and $z|d$
since $d$ is the gcd of a subset of the degrees.
If $d\in \md$ does not divide $\ell$, then $d|(\ell-1)r$ and $d$ is a gcd 
of a subset of $\set{r,\ldots,(\ell-1)r}$ so that $r|d$.
We have shown that
$\md\subseteq T_1\cup T_2$, where
\begin{equation}\label{eq:t1t2}
T_1=\set{d:z\vert d\vert\ell}\quad\hbox{and}\quad 
T_2=\set{kr:k|\ell-1}.
\end{equation}
Note that these sets intersect iff $\gcd(r,\ell)=r$, iff $r|\ell$.
To show the inclusion is an equality, we have to show
that $G(d)$ are distinct for distinct $d\in T_1\cup T_2$.
If $d\in T_1$, we have
$\gcd(r,d)=\gcd(r,\ell)=z$, so that by the proof of \cite[5.2]{lespr99},
$G(d)=G(dr/z,r,\ell z/d).$

On the other hand, if $d=kr$ where $k|\ell-1$ but $d\neq z$, we have
by the proof of \cite[Prop. 5.2]{lespr99} that 
$G(kr)=G(kr,1,(\ell-1)/k)$.
These groups have distinct parameters for each $d\in T_1\cup T_2$,
so
$\md=T_1 \cup T_2$ as claimed.
\end{proof}

\subsection{Proof of Theorem~\ref{th:cases} for $G(r,p,\ell)$}
It remains to show that $c(G)=0$ for irreducible $G=G(r,p,\ell)$,
except for the parameters $(de\ell,e\ell,\ell)\neq(3,3,3)$, and
for $(de\ell,e\ell,2\ell)$, where $e\ell$ is odd.  For these
exceptions, we show $c(G)=(-1)^{\ell-1}$.
It will be convenient to let
$S(m,k)=1/m\sum_{d|m}\mu(d)(-1)^{md/k-1}$, where $\mu$ is the 
(number-theoretic) M\"obius function on $\nat$.
\begin{lemma}
$$
S(m,k)=\begin{cases} (-1)^{km-1}, & \text{$m=1$, or $m=2$  and $k$ odd}\\
0,&\text{otherwise.}
\end{cases}
$$
\end{lemma}

\begin{proof}
Recall that $\sum_{d|m}\mu(d)=\delta_{m,1}$.
Consider the case where $k$ is odd and $m$ even.  Write $m=2^st$ for $t$ odd.
Then
\begin{eqnarray*}
S(m,k) &=& \frac{1}{m}\left(-\sum_{d|m,\text{$m/d$ even}}\mu(d)+
\sum_{d|m,\text{$m/d$ odd}}\mu(d)\right)\\
&=& \frac{1}{m}\left(-\sum_{d|m/2}\mu(d) -\sum_{d|m,2^s|d}\mu(2^s)
\mu(\frac{d}{2^s})\right)\\
&=& -\frac{2}{m}\delta_{m,2}=-\delta_{m,2}.
\end{eqnarray*}
The remaining cases are similar.
\end{proof}

At the same time, we will calculate $\chi_\Gamma$.
\begin{proposition}
Let $\Gamma=G\times C_m$ and $m=\deg(\delta_G)$.  Then
\begin{enumerate}
\item[(a)] For $G=G(r,p,{\ell})$ and $q=r/p>1$, 
$$
\chi_{\Gamma}=\begin{cases} 
I_{{\ell}q}, & \text{$\ell$ odd;}\\
-I_{{\ell}q}, & \text{$\ell$ and $p$ even;}\\
I_{{\ell}q}-I_{q\ell/2}, & \text{$\ell$ even, $p$ odd.}\\
\end{cases}
$$
\item[(b)] For $G=G(r,r,{\ell}),$ and ${\ell}>1$,
$$
\chi_{\Gamma}=
\begin{cases} 
I_{({\ell}-1)r} - I_{r({\ell}-1)/2} + I_{\ell}, & \text{$\ell$ odd;}\\
I_{({\ell}-1)r}-I_{\ell}, &\text{$r$, $\ell$ even;}\\
I_{({\ell}-1)r}+I_{\ell} -I_{\ell/2} &\text{$\ell$ even, $r$ odd.}\\
\end{cases}
$$
\end{enumerate}
\end{proposition}
\begin{proof}
We will handle $C_r$ as the special case $G(r,1,1)$ of (a)
and $S_{\ell}$ as the special case $G(1,1,{\ell})$ of (b).
Recall that $n=\rk(G(r,p,\ell))=\ell$ unless $r=p=1$ when 
$n=\rk(G(1,1,\ell))=\ell-1$ since $G(1,1,\ell)=S_\ell$ acts irreducibly
on an $\ell-1$ dimensional space.
\begin{enumerate}
\item[(a)] For $G=G(r,p,{\ell})$, $p<r$,
we have
$$
\md=\set{kq:t|k|{\ell}}
$$
where $t=\gcd(p,{\ell})$ and for $t|k|{\ell}$,
$G(kq)=G(kr/t,p,\ell t/k)$.
So 
\begin{eqnarray*}
c(G)&=&a_{tq} = tq\sum_{tq|kq}\mu(tq,kq)\frac{u(kq)}{i(kq)}\\
&=&\frac{tq}{{\ell}q}\sum_{t|k|{\ell}}\mu(k/t)(-1)^{\ell t/k-1}\\
&=&S(\ell/t,t)\\
&=&\begin{cases} 
(-1)^{\ell-1} & \text{$\ell=t$ or $\ell=2t$, $t$ odd ($\ell|p$ or 
$\frac{\ell}2|p$, $p$ odd.)}\\
0 & \text{otherwise.}\\
\end{cases}
\end{eqnarray*}

For $G(kq)=G(kr/t,p,\ell t/k)$ with $t|k|{\ell}$ we have
$$c(G(kq))=a_{kq}=\begin{cases}
(-1)^{\ell t/k-1}, &\text{$k=\ell$ or $k=\ell/2,$ $p$ odd}\\
0, &\text{otherwise.}
\end{cases}
$$
This shows that 
$$\chi_{\Gamma}=\begin{cases} 
I_{{\ell}q}, & \text{$\ell$ odd;}\\
-I_{{\ell}q}, & \text{$p$, ${\ell}$ even;}\\
I_{{\ell}q}-I_{q\ell/2}, & \text{$\ell$ even, $p$ odd,}\\
\end{cases}
$$
as required.

\item[(b)]
For $G=G(r,r,{\ell})$, $\ell>1$, by Proposition~\ref{prop:regnos},
$\md=T_1\cup T_2$, defined in \eqref{eq:t1t2}.
For $d\in T_1$, we have
$G(d)=G(dr/z,r,\ell z/d)$, and for $d\in T_2$, 
$G(d)=G(kr,1,(\ell-1)/k)$.
Since $z=\gcd(r,\ell)$, we have $z\in T_1$, and 
$z\in T_2$ if and only if $z=r$.  
\begin{eqnarray*}
c(G)&=&a_z=z\sum_{d\in \md}\mu(z,d)\frac{u(d)}{i(d)}\\
&=& z\frac{u(z)}{i(z)} + z\sum_{d\in T_1-\set{z}}\mu(z,d)\frac{u(d)}{i(d)}
+ z\sum_{kr\in T_2-\set{z}}\mu(z,kr)\frac{u(kr)}{i(kr)}\\
&=&
(-1)^{{\ell}-1}\frac{z\left(\left({\ell}-1\right)r-{\ell}\right)}{{\ell}
({\ell}-1)r} + 
\frac{z}{{\ell}}\sum_{1\ne\frac{d}{z}|\frac{{\ell}}{z}}\mu(\frac{d}{z})
(-1)^{\frac{{\ell}/z}{d/z}z-1}+\\
&&+ \frac{z}{({\ell}-1)r}
\sum_{k|{\ell}-1, kr\ne z}\mu(z,kr)(-1)^{(\ell-1)/k-1}\\
&=& (-1)^{{\ell}-1}\frac{z}{({\ell}-1)r} + S(\ell/z,z) +
\frac{z}{({\ell}-1)r}\sum_{k|{\ell}-1, kr\ne z}\mu(z,kr)(-1)^{(\ell-1)/k-1}
\end{eqnarray*}
Note that if $z\ne r$ then 
$\mu(z,r)=-1$ and $\mu(z,kr)=0$ for all $1\ne k|{\ell}-1$
whereas if $z=r$ then $\mu(z,kr)=\mu(kr/z)=\mu(k)$.
So we have
$$c(G)= S(\ell/z,z)+\delta_{zr}S({\ell}-1,1)$$
This means that $c(G)=(-1)^{n-1}$ where $n=\rk(G)$ iff
$$
\begin{cases} 
{\ell}|r, &(r,{\ell})\ne (2,2),(3,3)\\
\frac{{\ell}}{2}|r &r\mbox{ odd}\\
(r,{\ell})=(1,3)
\end{cases}
$$
Otherwise $c(G)=0$.
Note that $c(G(2,2,2))=0$ agrees with the statement since $G(2,2,2)$
is reducible. Also observe that if $G=G(r,p,{\ell})$ is a rank 2 irreducible,
we obtain $c(G)=-1$ as was predicted in Section~\ref{sec:rank2}.
This includes the case of the rank 2 irreducible $G(1,1,3)=S_3$. 

To compute $\chi_\Gamma$, it remains only to find $a_d$, for $d\in\md$.
If $d\in T_1$, by part (a), 
$$a_d=c(G(\frac{dr}{z},r,\frac{{\ell}z}{d}))=
\begin{cases}
(-1)^{\frac{{\ell}z}{d}-1} &\text{$d=\ell$, or $d={\ell}/2$, $r$ odd;}\\
0 &\text{otherwise.}
\end{cases}
$$
For $d\in T_2$, 
$$a_d=c(G(kr,1,\frac{{\ell}-1}{k}))
=\begin{cases} 
1 &\text{$d=({\ell}-1)r$;}\\
-1 &\text{$d=r(\ell-1)/2$, $\ell$ odd;}\\
0 &\text{otherwise.}
\end{cases}
$$
The expression for $\chi_{\Gamma}$ in (b) follows. 
\end{enumerate}
\end{proof}

\begin{figure}
$$
\begin{array}{|c|l||c|l|} \hline
\hbox{\#}&\chi_\Gamma & \hbox{\#}&\chi_\Gamma\\ \hline
3 & I_r             & 13 & I_{12}-I_4 \\        
4 & I_6+I_4-I_2     & 14 & I_{24}-I_6 \\        
5 & I_{12}-I_6      & 15 & -I_{12} \\          
6 & I_{12}-I_4      & 16 & I_{30}+I_{20}-I_{10} \\
7 &  -I_{12}        & 17 & I_{60}-I_{20} \\    
8 & I_{12}+I_8-I_4  & 18 & I_{60}-I_{30} \\    
9 & I_{24}-I_8      & 19 & -I_{60} \\          
10 & I_{24}-I_{12}  & 20 & I_{30}+I_{12}-I_6 \\
11 & -I_{24}        & 21 & I_{60}-I_{12} \\     
12 &I_8+I_6-I_2     & 22 & I_{20}+I_{12}-I_4 \\ \hline
\end{array}
$$
\caption{$\chi_\Gamma$ for $G_3$ and rank-$2$ exceptionals}
\label{figa}
\end{figure}
\begin{figure}
$$
\begin{array}{|c|l||c|l|} \hline
\hbox{\#}&\chi_\Gamma & \hbox{\#}&\chi_\Gamma\\ \hline
23 & I_{10}+I_6     & 28 & I_{12}+I_8-I_6-I_4 \\                     
24 & I_{14}+I_6     & 29 & I_{20}-I_4 \\                             
25 & I_{12}+I_9-I_6 & 30 & I_{30}+I_{20}+I_{12}-I_{10}-I_6-I_4 \\
26 & I_{18}	    & 31 & I_{24}+I_{20}-I_{12}-I_8 \\
27 & I_{30}	    & 32 & I_{30}+I_{24}-I_{12} \\ \hline
\end{array}
$$
\caption{$\chi_\Gamma$ for exceptionals of rank $3$ and $4$}
\label{figb}
\end{figure}
\begin{figure}
$$
\begin{array}{|c|c|l|} \hline
\hbox{\#}&\hbox{rank}&\chi_\Gamma \\ \hline
33 & 5 & I_{18}+I_{10} \\
34 & 6 & I_{42}-I_6 \\
35 & 6 & I_{12}+I_8-I_6-I_4 \\
36 & 7 & I_{18}+I_{14} \\
37 & 8 & I_{30}+I_{24}+I_{20}-I_{12}-I_{10}-I_8\\ \hline
\end{array}
$$
\caption{$\chi_\Gamma$ for exceptionals of rank $\geq5$}
\label{figc}
\end{figure}
\section{orbifold Euler characteristics}
The orbifold Euler characteristic of a space $X$ under the action of 
a group $G$ is defined to be
$$
\sum_{[g]}\chi(X^g/C_G(g)),
$$
where the sum is taken over all conjugacy classes of $G$.

The orbifold Euler characteristic of the Milnor fibre $F$ under the
action of a $G\times C_m$ can be expressed in terms of the integers
$\set{a_d}$ from \eqref{eq:defofa}, and we include its calculation here
as an example.

\begin{lemma}\label{lem:orbi}
For any reflection group $G$,
$$
\chi(U/G)=\sum_{d\in\md}a_d.
$$
\end{lemma}
\begin{proof}
We have a disjoint union $U=\bigcup_{d\in\md}G\cdot U(d)^\circ$.
  The claim follows
by Proposition~\ref{prop:freely}, Lemma~\ref{lem:ecs}, 
and the definition $a_d=\chi(U(d)^\circ)/G(d)$.
\end{proof}

\begin{theorem}
Let $G$ be a reflection group, $f$ an unreduced
discriminant polynomial of degree $m$,
and $F$ its Milnor fibre.
The orbifold Euler characteristic of $F$ with respect to $\Gamma=G\times C_m$
equals
$$
\sum_{d\in\md}d a_d,
$$
in the notation of Section~\ref{sec:ec}.
\end{theorem}
\begin{proof}
For a regular pair $(g,\zeta)$, 
we have $F^{(g,\zeta)}/C(g,\zeta)=U^g/C_G(g)$.  Recall that if 
$(g,\zeta)$ is not a regular pair, then the set of fixed points is 
empty, so we need only consider a sum over conjugacy classes of 
regular pairs in $G\times C_m$.

Using Lemma~\ref{lem:id}, there are $\phi(d)$ conjugacy classes of
regular pairs of order $d$, for each regular number $d$.
So we have
\begin{eqnarray*}
\sum_{[(g,\zeta)]}\chi(F^{(g,\zeta)}/C(g,\zeta))&=&
\sum_{d\in\Reg}\phi(d) \chi(U^g/C_G(g))\\
&=&\sum_{d\in\Reg}\phi(d)a_{\roundup{d}}\\
&=& \sum_{d\in\md}\sum_{d'|d}\phi(d')a_d\\
&=& \sum_{d\in\md}d a_d,
\end{eqnarray*}
where the second equality follows from Lemma~\ref{lem:orbi} together
with Lemma~\ref{lem:czrs}.
\end{proof}
\begin{rem}
By way of comparison, the ordinary or orbifold Euler characteristics
of $F/\Gamma$ are equal to the image of $\chi_\Gamma$ under 
homomorphisms from the character ring of $\Gamma$ to $\Z$ that take
$I_d$ to $1$, or to $d$, respectively.
\end{rem}
\section{Concluding remarks}
This investigation leaves the obvious open question of 
whether Theorem~\ref{th:cases} could be proven in a more conceptual way.
Our proof depends on knowing (co)degrees and regular numbers for each
group, which are not reflected in the simplicity of the statement.  

We also note that $\Gamma=G\times C_m$ is not the most general group for
which these calculations make sense.  In general one should replace 
$G$ by $N(G)$, the normalizer of $G$ in $U(V)$, and $C_m$ by 
$\Lambda= C_m\rtimes{\rm Gal}(K_m/K)$, where $K$ is the splitting field for $N(G)$
and $K_m$ is the extension of $K$ containing all $m$th roots of unity.
Note that ${\rm Gal}(K_m/K)$ is a finite group of order dividing
$\phi(m)$ which acts on $\C$ by inflation to ${\rm Gal}(\C/K)$,
since $K_m/K$ is a Galois extension. This action can be extended
to a diagonal action on  $\C^n$ which stabilizes $F=\delta^{-1}(1)$ since the
coefficients of $\delta$ lie in $K$. Note that ${\rm Gal}(K_m/K)\cap C_m=1$
and that ${\rm Gal(K_m/K)}$ normalizes $C_m$.
The actions of ${\rm Gal}(K_m/K)$ and
$N(G)$ on $F$ commute by construction.  
It may be interesting to examine the $\Lambda$-module structure of the
equivariant Euler characteristic $\chi_\Lambda$.
It is probable that an answer would 
involve Springer's twisted regular numbers.   
   
\providecommand{\bysame}{\leavevmode\hbox to3em{\hrulefill}\thinspace}


\begin{thebibliography}{10}

\bibitem{bjo94}
A.~Bj{\"o}rner, \emph{Subspace arrangements}, First European Congress of
  Mathematics, Vol.\ I (Paris, 1992), Birkh\"auser, Basel, 1994, pp.~321--370.

\bibitem{bmr98}
M.~Brou{\'e}, G. Malle, and R. Rouquier, \emph{Complex
  reflection groups, braid groups, {H}ecke algebras}, J. Reine Angew. Math.
  \textbf{500} (1998), 127--190.

\bibitem{kbrown82}
K.~S. Brown, \emph{Complete {E}uler characteristics and fixed-point
  theory}, J. Pure Appl. Algebra \textbf{24} (1982), no.~2, 103--121.

\bibitem{co76}
A.~M. Cohen, \emph{Finite complex reflection groups}, Ann. Sci. \'Ecole
  Norm. Sup. (4) \textbf{9} (1976), no.~3, 379--436.

\bibitem{crapo}
H.~Crapo, \emph{A higher invariant for matroids}, J. Comb. Th. \textbf{2}
  (1967), 406--417.

\bibitem{CRvolI}
C.~W. Curtis and I. Reiner, \emph{Methods of representation theory.
  {V}ol. {I}}, John Wiley \& Sons Inc., New York, 1990, With applications to
  finite groups and orders, Reprint of the 1981 original, A Wiley-Interscience
  Publication.

\bibitem{cpss99}
C.~De~Concini, C. Procesi, M. Salvetti, and F. Stumbo,
  \emph{Arithmetic properties of the cohomology of {A}rtin groups}, Ann. Scuola
  Norm. Sup. Pisa Cl. Sci. (4) \textbf{28} (1999), no.~4, 695--717.

\bibitem{delo95}
J.~Denef and F.~Loeser, \emph{Regular elements and monodromy of discriminants
  of finite reflection groups}, Indag. Math. (N.S.) \textbf{6} (1995), no.~2,
  129--143.

\bibitem{d97}
G.~Denham, \emph{{H}anlon and {S}tanley's conjecture and the {M}ilnor fibre of
  a braid arrangement}, J. Algebraic Combin. \textbf{11} (2000), 227--240.

\bibitem{de01}
G. Denham, \emph{The combinatorial {L}aplacian of the {T}utte complex}, J.
  Algebra \textbf{242} (2001), no.~1, 160--175.

\bibitem{dur79}
A.~H. Durfee, \emph{Fifteen characterizations of rational double points and
  simple critical points}, Enseign. Math. (2) \textbf{25} (1979), no.~1-2,
  131--163.

\bibitem{klein}
F. Klein, \emph{Lectures on the icosahedron and the solution of equations of
  the fifth degree}, revised ed., Dover Publications Inc., New York, N.Y.,
  1956, Translated into English by George Gavin Morrice.

\bibitem{leh95}
G.~I. Lehrer, \emph{Poincar\'e polynomials for unitary reflection groups},
  Invent. Math. \textbf{120} (1995), no.~3, 411--425.

\bibitem{lespr99}
G.~I. Lehrer and T.~A. Springer, \emph{Reflection subquotients of unitary
  reflection groups}, Canad. J. Math. \textbf{51} (1999), no.~6, 1175--1193,
  Dedicated to H. S. M. Coxeter on the occasion of his 90th birthday.

\bibitem{milnorbook}
J.~W. Milnor, \emph{Singular points of complex hypersurfaces}, Princeton
  University Press, 1968.

\bibitem{ot}
P.~Orlik and H.~Terao, \emph{Arrangements of hyperplanes}, Grundlehren der
  Mathematischen Wissenschaften, no. 300, Springer-Verlag, 1992.

\bibitem{os80b}
P.~Orlik and L.~Solomon, \emph{Unitary reflection groups and cohomology},
  Invent. Math. \textbf{59} (1980), no.~1, 77--94.

\bibitem{os84}
\bysame, \emph{Arrangements defined by unitary reflection groups}, Math. Ann.
  \textbf{261} (1982), no.~3, 339--357.

\bibitem{spr74}
T.~A. Springer, \emph{Regular elements of finite reflection groups}, Invent.
  Math. \textbf{25} (1974), 159--198.

\end{thebibliography}
\end{document}